\documentclass[12pt]{article}

\begin{document}
\bibliographystyle{plain}

\thispagestyle{empty}
\setcounter{page}{0}

{\large Guszt\'av MORVAI and  Benjamin WEISS: }

\vspace {2cm}

{\Large Order Estimation of Markov Chains}

\vspace {2cm}

{\large IEEE Trans. Inform. Theory  51  (2005),  no. 4, 1496--1497.}

\vspace {2cm}

\begin{abstract}
We describe estimators $\chi_n(X_0,X_1,\dots,X_n)$, which when applied to an unknown stationary 
process taking values from  a countable alphabet ${\cal X}$, converge almost surely to $k$ 
in case the process is 
a $k$-th order Markov chain and to infinity otherwise. 
\end{abstract}



\pagebreak

\section{Introduction}

When faced with an unknown stationary and ergodic stochastic process $X_1,X_2,\dots, X_n,\dots$ one 
may try to determine various properties of this process from the successive observations up to time $n$. 
For example, one might try to estimate the entropy of the process. Several schemes of the form 
$g_n(X_1,\dots,X_n)$ are known which will converge almost surely to the entropy of the process $\{X_n\}$
cf. Bailey \cite{Bailey76}, Csisz\'ar and Shields \cite{CsSh00}, Csisz\'ar \cite{Cs02}, Ornstein and Weiss
\cite{OW93}, \cite{OW90}, \cite{OW02}, Kontoyiannis,  Algoet,  Suhov and  Wyner \cite{KASW98}  and Ziv \cite{Ziv78}.
However, if one just wants to determine whether or not the process has positive entropy 
(often associated with the popular notion of chaos) then there is no sequence of
two valued functions $e_n(X_1,\dots,X_n)\in \{ZERO,POSITIVE\}$ with the property that almost surely, $e_n$ 
stabilize at $ZERO$  
for all zero entropy processes and at $POSITIVE$ for all positive entropy processes. 
(While this result does not appear explicitly in Ornstein amd Weiss \cite{OW90}, it can be readily established using a 
very simple  
variant 
of the construction given there in $\S$ $4$.)

A similar situation obtains in testing for membership in the class of  $k$-th order Markov chains. 
One can estimate the order of a Markov chain by e.g the method of 
Csisz\'ar and Shields \cite{CsSh00} or  Csisz\'ar \cite{Cs02}. 
They show that the minimum description length Markov estimator will converge almost surely to the correct order if 
the alphabet size is bounded a priori. Without this assumption they show that this is no longer true. To accomplish 
their goals they study the large scale typicality of Markov sample paths. A further negative result is that of 
Bailey \cite{Bailey76} who showed that no two valued test exists for testing mixing Markov vs. not mixing Markov. 

We will present a more direct estimator for the order of a Markov chain which also uses the fact that there are 
universal rates for the convergence of empirical $k$-block distributions in this class. Our approach enables us to 
dispense  with  the assumption 
that the alphabet size is bounded, indeed it may even be infinite, as long as  there is a 
finite memory. In addition we will show 
that if the process is not a Markov chain then the estimate for the order will tend to infinity. 
This is in complete analogy with the entropy estimation that we mentioned earlier.

\section{The Order Estimator}

\smallskip
\noindent
Let $\{X_n\}_{n=-\infty}^{\infty}$ be a stationary and ergodic time series taking values from a 
discrete (finite or countably infinite)  alphabet 
${\cal X}$. (Note that all stationary time series $\{X_n\}_{n=0}^{\infty}$ 
can be thought to be a 
two sided time series, that is, $\{X_n\}_{n=-\infty}^{\infty}$. )  
For notational convenience, let $X_m^n=(X_m,\dots,X_n)$,
where $m\le n$. Note that if $m>n$ then $X_m^n$ is the empty string. 

\noindent
Let 
  $p(x_{-k}^0)$ and $p(y|x^0_{-k})$ denote the  distribution $P(X_{-k}^{0}=x_{-k}^{0})$ and  
the conditional distribution $P(X_1=y|X^0_{-k}=x^0_{-k})$, respectively.

\noindent
A discrete alphabet stationary time series is said to be a Markov chain if 
for some $K\ge 0$, for all $y\in {\cal X}$,  $i\ge 1$ and $z^0_{-K-i+1}\in {\cal X}^{K+i}$, 
if $p( z^0_{-K-i+1})>0$ then 
$$
 p(y| z^0_{-K+1})= p(y| z^0_{-K-i+1}).
$$
The order of a Markov chain is the smallest such $K$. 

\bigskip
\noindent
In order to estimate the order  we need to define some explicit statistics. 

\smallskip
\noindent
For $k\ge 0$ let ${\cal S}_{k}$ denote the support  of the distribution of $X^0_{-k}$ as    
$${\cal S}_{k}=
\{x^0_{-k}\in {\cal X}^{k+1}:
 p(x^0_{-k})>0\}.
$$ 

\smallskip
\noindent
Define
$$
\Delta_k=\sup_{1\le i} \sup_{(z^0_{-k-i+1},x)\in {\cal S}_{k+i}}
\left| p(x| z^0_{-k+1})- p(x|z^{0}_{-k-i+1})\right|.
$$

\noindent
We will divide the data segment $X_0^n$ into two parts: $X_0^{\lceil{n\over 2}\rceil-1}$ and 
$X_{\lceil {n\over 2} \rceil}^n$. Let ${\cal S}_{n,k}^{(1)}$ denote the set of strings with length $k+1$ 
which appear at all in 
$X_0^{\lceil{n\over 2}\rceil-1}$. That is,
$$
{\cal S}_{n,k}^{(1)}= \{x^0_{-k}\in {\cal X}^{k+1}: 
 \exists k\le t \le \lceil{n\over 2}\rceil-1 : X^t_{t-k}=x^0_{-k}\}.
$$ 

\noindent
For a fixed  $0<\gamma<1$  let ${\cal S}_{n,k}^{(2)}$ denote the set of strings with length $k+1$  which appear  
more than $n^{1-\gamma}$ times in $X_{\lceil {n\over 2} \rceil}^n$. That is, 
$${\cal S}_{n,k}^{(2)}=\{x^0_{-k}\in {\cal X}^{k+1}: 
\#\{\lceil {n\over 2} \rceil+k\le t\le n: X^t_{t-k}=x^0_{-k}\} > n^{1-\gamma}\}.$$
Let 
$$
{\cal S}_{k}^n={\cal S}_{n,k}^{(1)}\bigcap {\cal S}_{n,k}^{(2)}.
$$
For notational convenience, let $C(x|z^0_{-k+1}: [n_1,n_2])$ denote the empirical conditional probability 
of $X_1=x$ given $X^0_{-k+1}=z^0_{-k+1}$ from the samples $(X_{n_1},\dots,X_{n_2})$, that is,   
$$
C(x|z^0_{-k+1}: [n_1,n_2])=
{ \#\{ n_1+k\le t\le n_2: X^t_{t-k}=(z^0_{-k+1},x)\}\over 
\#\{n_1+k-1\le t\le n_2-1: X^t_{t-k+1}=z^0_{-k+1}\}}
$$
where $0/0$ is defined as $0$.

\noindent
 We define the empirical version of $\Delta_k$ as follows:  
$$
 {\hat \Delta}^n_k=\max_{1\le i \le n}
\max_{(z^{0}_{-k-i+1},x)\in {\cal S}^n_{k+i} } 
 \left|
 C(x|z^0_{-k+1}: [\lceil {n\over 2} \rceil,n])-
  C(x|z^0_{-k-i+1}: [\lceil {n\over 2} \rceil,n])
\right|.  
$$

\smallskip
\noindent
Observe, that by ergodicity, for any fixed $k$, 
\begin{equation}\label{Deltatozero}
\liminf_{n\to\infty}{\hat \Delta}^n_k\ge \Delta_k \ \ \mbox{almost surely.}
\end{equation}

\bigskip
\noindent
We define an estimate $\chi_n$ for the order   from samples 
$X_0^n$ as follows.  
Let $0< \beta <{1-\gamma \over 2}$ be  arbitrary. Set $\chi_0=0$, and for $n\ge 1$ 
let  $\chi_n$ be the smallest 
$0\le k_n< n $ such that 
${\hat \Delta }^n_{k_n}\le n^{-\beta}$.

\bigskip
\noindent
{\bf  THEOREM. }
{\it If the stationary and ergodic time series $\{X_n\}$ taking values from a discrete alphabet 
happens to be a Markov chain with any finite order  then 
$\chi_n$ equals to the order eventually almost surely, 
and if it is not Markov with any finite order then  
$\chi_n\to\infty$ 
almost surely. 
}

\bigskip
\noindent
{\bf Application:} 
Let $M> 0$ be arbitrary. 
The goal is to  decide if the discrete alphabet stationary and ergodic 
time series is a Markov chain with order less than   $M$ or not. One may use $\chi_n$ 
and say YES if $\chi_n< M$ and say NO otherwise. By the  Theorem, 
eventually, the answer will be correct. 
 
\section{Proof of the Theorem}

\bigskip
\noindent
{\bf Proof:} 
If the process is a Markov chain, it is immediate that for all 
$k$ greater than or equal  the order, $\Delta_k=0$. For $k$ less than   the order $ \Delta_k>0$.   
If the process is not a Markov chain with any finite order  then $\Delta_k>0$ 
for all $k$.
Thus by (\ref{Deltatozero}) if the process is not Markov  then $\chi_n\to\infty$ and if 
 it is Markov  then 
$\chi_n$ is greater or equal  the order eventually almost surely. 
We have to show that $\chi_n$ is less or equal  the order  eventually almost surely 
provided that the process is a Markov chain. 
 
 \noindent
Assume that the process is a Markov chain with order $k$. Let $n\ge k$. 
We will estimate  the probability of the undesirable  event as follows:
\begin{eqnarray*}
\lefteqn{
P({\hat \Delta}^n_k> n^{-\beta}|X_0^{\lceil {n\over 2} \rceil})\le }\\
&&\sum_{i=1}^{n}
  P(\max_{(z^0_{-k-i+1},x)\in {\cal S}^n_{k+i} }  \left| 
 C(x|z^0_{-k+1}: [\lceil {n\over 2} \rceil,n])-
 C(x|z^0_{-k-i+1}: [\lceil {n\over 2} \rceil,n])
\right|> n^{-\beta}|X_0^{\lceil {n\over 2} \rceil}).
 \end{eqnarray*}
 We can estimate each  probability in the sum as the sum of two terms:

\begin{eqnarray*}
\lefteqn{ P(\max_{(z^{0}_{-k-i+1},x)\in {\cal S}^n_{k+i} }\left| 
C(x|z^0_{-k+1}: [\lceil {n\over 2} \rceil,n])-
 C(x|z^0_{-k-i+1}: [\lceil {n\over 2} \rceil,n])\right|
 > n^{-\beta} |X_0^{\lceil {n\over 2} \rceil}) }\\
&\le& 
 P(\max_{(z^{0}_{-k-i+1},x)\in {\cal S}^n_{k+i} }\left| 
C(x|z^0_{-k+1}: [\lceil {n\over 2} \rceil,n])-
 p(x|z^0_{-k+1})\right|
 > 0.5 n^{-\beta} |X_0^{\lceil {n\over 2} \rceil})\\
&+& 
P(\max_{(z^{0}_{-k-i+1},x)\in {\cal S}^n_{k+i} }\left|
p(x|z^0_{-k+1})-
C(x|z^0_{-k-i+1}: [\lceil {n\over 2} \rceil,n])\right|
 >0.5 n^{-\beta}|X_0^{\lceil {n\over 2} \rceil}).
\end{eqnarray*}
We overestimate these probabilities. 
For any $m\ge 0$ and $x^0_{-m}$ define  $\sigma^m_i(x^0_{-m})$ as the time of the $i$-th ocurrence of the string $x^0_{-m}$ in 
 the data segment $X_{\lceil {n\over 2} \rceil}^n$, that is, let $\sigma^m_{0}(x^0_{-m})=\lceil {n\over 2} \rceil+m-1$ 
 and for $i\ge 1$ define
 $$
 \sigma^m_i(x^0_{-m})=\min\{t>\sigma^m_{i-1}(x^0_{-m}) : X^t_{t-m}=x^0_{-m}\}.
 $$
 Now  
 \begin{eqnarray*}
\lefteqn{ 
P(\max_{(z^{0}_{-k-i+1},x)\in {\cal S}^n_{k+i} }\left| 
C(x|z^0_{-k+1}: [\lceil {n\over 2} \rceil,n])-
 C(x|z^0_{-k-i+1}: [\lceil {n\over 2} \rceil,n])\right|
 > n^{-\beta} |X_0^{\lceil {n\over 2} \rceil})}\\
&\le& 
 P(\max_{(z^0_{-k+1},x)\in {\cal S}^{(1)}_{n,k} }\sup_{j> n^{1-\gamma}}
 \\
&& \left| {1\over j} \sum_{r=1}^j 1_{\{X_{\sigma_r^{k-1}(z^0_{-k+1})}=x\}} -p(x|z^0_{-k+1})\right| 
  >0.5 n^{-\beta}|X_0^{\lceil {n\over 2} \rceil})\\
&+&
  P(\max_{(z^{0}_{-k-i+1},x)\in {\cal S}^{(1)}_{n,k+i} }\sup_{j> n^{1-\gamma}} \\
&& \left| {1\over j} \sum_{r=1}^j 1_{\{X_{\sigma_r^{k+i-1}(z^{0}_{-k-i+1})}=x\}} -
p(x|z^0_{-k+1})\right| 
  >0.5 n^{-\beta}|X_0^{\lceil {n\over 2} \rceil})
\end{eqnarray*}
Since both ${\cal S}^{(1)}_{n,k}$ and ${\cal S}^{(1)}_{n,k+i}$ depend solely on $X_0^{\lceil {n\over 2} \rceil}$ we get 
\begin{eqnarray*}
\lefteqn{P(\max_{(z^{0}_{-k-i+1},x)\in {\cal S}^n_{k+i} }\left| 
C(x|z^0_{-k+1}: [\lceil {n\over 2} \rceil,n])-
 C(x|z^0_{-k-i+1}: [\lceil {n\over 2} \rceil,n])\right|
 > n^{-\beta} |X_0^{\lceil {n\over 2} \rceil})}\\
&\le& \sum_{(z^0_{-k+1},x)\in {\cal S}^{(1)}_{n,k} }
\sum_{j= \lceil n^{1-\gamma}\rceil}^{\infty}
   P(
\left| {1\over j} \sum_{r=1}^j 1_{\{X_{\sigma_r^{k-1}(z^0_{-k+1})}=x\}} -p(x|z^0_{-k+1})\right|\\
&&  >0.5 n^{-\beta}|X_0^{\lceil {n\over 2} \rceil})\\
  &+&\sum_{(z^{0}_{-k-i+1},x)\in {\cal S}^{(1)}_{n,k+i} }
  \sum_{j= \lceil n^{1-\gamma}\rceil}^{\infty}
   P( \left| {1\over j} \sum_{r=1}^j 1_{\{X_{\sigma_r^{k+i-1}(z^{0}_{-k-i+1})}=x\}} \right.\\
  &&  \left. -p(x|z^0_{-k+1})
   \right|
  >0.5 n^{-\beta}|X_0^{\lceil {n\over 2} \rceil}).
  \end{eqnarray*}
   Each of these represents the deviation of an empirical count from its mean. 
The variables in question are independent since whenever the block 
 $z^0_{-k+1}$ occurs the next term is chosen using the same distribution $p(x|z^0_{-k+1})$.
 Thus by Hoeffding's inequality (cf. Hoeffding \cite{Hoeffding63} or Theorem 8.1 of Devroye et. al. \cite{DGYL96}) for 
 sums of bounded independent random variables and since the cardinality of both 
 ${\cal S}^{(1)}_{n,k}$ and ${\cal S}^{(1)}_{n,k+i}$ is not greater than $ (n+2)/2$, 
we have 
\begin{eqnarray*}
 &&P(\max_{(z^{0}_{-k-i+1},x)\in {\cal S}^n_{k+i} }\left| 
C(x|z^0_{-k+1}: [\lceil {n\over 2} \rceil,n])-
 C(x|z^0_{-k-i+1}: [\lceil {n\over 2} \rceil,n])\right|
 > n^{-\beta} |X_0^{\lceil {n\over 2} \rceil})\\
&&\le  2 {n+2\over 2} \sum_{j= \lceil n^{1-\gamma}\rceil}^{\infty} 2 e^{-2n^{-2\beta}j}.
\end{eqnarray*}
Thus 
$$
P({\hat \Delta}^n_k> n^{-\beta}|X_0^{\lceil {n\over 2} \rceil})\le  n (n+2)   4 e^{-2n^{-2\beta+1-\gamma}}.
$$
Integrating both sides we get 
$$
P({\hat \Delta}^n_k> n^{-\beta})\le  n (n+2)   4 e^{-2n^{-2\beta+1-\gamma}}.
$$
The right hand side is summable provided $2\beta+\gamma<1$ and 
the Borel-Cantelli Lemma yields that  
$P({\hat \Delta}^n_k\le n^{-\beta}  \ \mbox{eventually})=1$.
Thus  
$\chi_n\le  k$ eventually almost surely provided the process is Markov with  order  $k$.
The proof of the Theorem is complete.



\begin{thebibliography}{00}



\bibitem{Bailey76}
D. H. Bailey,
{\it Sequential Schemes for Classifying and Predicting 
Ergodic Processes.} Ph. D. thesis, Stanford University, 1976.

\bibitem{CsSh00}
I. Csisz\'ar and P. Shields,
"The consistency of the BIC Markov order estimator,"
{\it Annals of Statistics.}, vol. 28, pp. 1601-1619, 2000. 

\bibitem{Cs02}
I. Csisz\'ar,
"Large-scale typicality of Markov sample paths and consistency of MDL order estimators ,"
{\it IEEE Transactions on Information Theory}, vol. 48, pp. 1616-1628, 2002. 


\bibitem{DGYL96}
L Devroye, L. Gy\"orfi, G. Lugosi,
{\it A Probabilistic Theory of Pattern Recognition.} Springer-Verlag, New York, 1996.




\bibitem{Hoeffding63}
W. Hoeffding,
"Probability inequalities for sums of bounded random variables ,"
{\it Journal of the American Statistical Association}, vol. 58, pp. 13-30, 1963. 

\bibitem{KASW98}
I. Kontoyiannis, P. Algoet, Yu.M. Suhov, A.J. Wyner, 
"Nonparametric entropy estimation for stationary processes and random fields, with application to English text,"
{\it IEEE Transactions on Information Theory}, vol. 44, pp. 1319--1327, 1998.

\bibitem{OW90}
D. S. Ornstein and B. Weiss,
"How sampling reveals a process,"
{\it The Annals of Probability}, vol. 18, pp. 905--930, 1990.


\bibitem{OW93}
D. S. Ornstein and B. Weiss,
"Entropy and data compression schemes,"
{\it IEEE Transactions on Information Theory}, vol. 39, pp. 78--83, 1993.


\bibitem{OW02}
D. S. Ornstein and B. Weiss,
"Entropy and recurrence rates for stationary random fields,"
{\it IEEE Transactions on Information Theory}, vol. 48, pp. 1699--1697, 2002.


\bibitem{Ziv78}
J. Ziv, " Coding theorems for individual sequences. 
{\it IEEE Transactions on Information Theory}, vol. 24, pp. 405--412, 1978.

\end{thebibliography}
\end{document}